\documentclass[reqno]{amsart}

\usepackage{amsmath,amsthm,amssymb,amsfonts,enumerate,graphicx,color,pstricks}
\usepackage[utf8]{inputenc}

\numberwithin{equation}{section} 
\theoremstyle{plain}
\newtheorem{theo+}           {Theorem}      [section]
\newtheorem{prop+}  [theo+]  {Proposition}
\newtheorem{coro+}  [theo+]  {Corollary}
\newtheorem{lemm+}  [theo+]  {Lemma}
\newtheorem{defi+}  [theo+]  {Definition}
\newtheorem{conj+}  [theo+]  {Conjecture}

\theoremstyle{definition}
\newtheorem{rema+}  [theo+]  {Remark}
\newtheorem{prob+}  [theo+]  {Problem}
\newtheorem{exam+}  [theo+]  {Example}

\newenvironment{corollary}{\begin{coro+}}{\end{coro+}}
\newenvironment{lemma}{\begin{lemm+}}{\end{lemm+}}

\newcommand{\ct}{\operatorname{CT}}

\makeatletter
\@namedef{subjclassname@2020}{\textup{2020} Mathematics Subject Classification}
\makeatother

\begin{document}

\baselineskip 18pt
\larger[2]
\title
[New proofs of the septic Rogers--Ramanujan identities] 
{New proofs of the septic Rogers--Ramanujan identities}
\author{Hjalmar Rosengren}
\address
{Department of Mathematical Sciences
\\ Chalmers University of Technology and University of Gothenburg\\SE-412~96 G\"oteborg, Sweden}
\email{hjalmar@chalmers.se}
\urladdr{http://www.math.chalmers.se/{\textasciitilde}hjalmar}
\thanks{Supported by the Swedish Science Research Council, project no.\ 2020-04221.}

\maketitle

\begin{abstract}
We give new proofs of the twelve Rogers--Ramanujan-type identities due to Rogers and Slater that are traditionally associated with the moduli $7$, $14$ and $28$. 
\end{abstract}

\section{Introduction}

The celebrated Rogers--Ramanujan identities 
\begin{subequations}\label{rri}
\begin{align}G(q)&=\sum_{k=0}^\infty\frac{q^{k^2}}{(q;q)_k}=\frac 1{(q,q^4;q^5)_\infty},\\
H(q)&=\sum_{k=0}^\infty\frac{q^{k(k+1)}}{(q;q)_k}=\frac 1{(q^2,q^3;q^5)_\infty},
\end{align}
\end{subequations}
were first proved by Rogers in  1894 \cite{r94} (see \S \ref{ps} for the notation). He obtained many similar results, such as
$$\sum_{k=0}^\infty\frac{q^{2k^2}}{(q^2;q^2)_k(-q;q)_{2k}}=\frac{(q^3,q^4,q^7;q^7)_\infty}{(q^2;q^2)_\infty}, $$
which is  related to modulus $7$ rather than $5$. In total,
Rogers found nine Rogers--Ramanujan-identities
 associated with moduli $7$ and $14$  \cite{r94,r17}, and  Slater's list \cite{slater}
contains in addition three modulus $28$ identities. Since the assignment of moduli 
 is partly a matter of convention,  we will refer to all twelve identities as \emph{septic}.

We recently gave a new proof of \eqref{rri} \cite{r22}, based on the identities
\begin{subequations}\label{rrr}
\begin{align}
G(q)&=\frac{(q^8;q^8)_\infty}{(q^2;q^2)_\infty}(G(q^{16})+qH(-q^4)),\\
H(q)&=\frac{(q^8;q^8)_\infty}{(q^2;q^2)_\infty}(G(-q^4)+q^3H(q^{16}).
\end{align}
\end{subequations}
 This can be viewed as a recursion for computing the power series coefficients of $G$ and $H$. To prove  \eqref{rri} it is therefore enough to verify that \eqref{rrr} holds for both the sum sides and the product sides. Surprisingly, this is essentially contained in  Rogers' 1894 memoir \cite{r94}, although he seems not to have realized that it immediately implies \eqref{rri}.  The proof given in \cite{r94} is much more complicated.

In the present work, we will apply the method of \cite{r22} 
to give new proofs of all 
 twelve septic Rogers--Ramanujan identities, which can be stated as
\begin{subequations}\label{ssa}
\begin{align}
\label{s33}\sum_{k=0}^\infty\frac{q^{2k^2}}{(q^4;q^4)_k(-q;q^2)_k}&=\frac{(q^3,q^4,q^7;q^7)_\infty}{(q^2;q^2)_\infty}, \\
\label{s32}\sum_{k=0}^\infty\frac{q^{2k(k+1)}}{(q^4;q^4)_k(-q;q^2)_k}&=\frac{(q^2,q^5,q^7;q^7)_\infty}{(q^2;q^2)_\infty},  \\
\label{s31}\sum_{k=0}^\infty\frac{q^{2k(k+1)}}{(q^4;q^4)_k(-q;q^2)_{k+1}}&=\frac{(q,q^6,q^7;q^7)_\infty}{(q^2;q^2)_\infty},  \\
\label{s61}\sum_{k=0}^\infty\frac{q^{k^2}}{(q;q)_k(q;q^2)_k}&=\frac{(q^6,q^8,q^{14};q^{14})_\infty}{(q;q)_\infty}, \\
\label{s60} 
\sum_{k=0}^\infty\frac{q^{k(k+1)}}{(q;q)_k(q;q^2)_{k+1}}&=\frac{(q^4,q^{10},q^{14};q^{14})_\infty}{(q;q)_\infty}, \\
\label{s59}\sum_{k=0}^\infty\frac{q^{k(k+2)}}{(q;q)_k(q;q^2)_{k+1}}&=\frac{(q^2,q^{12},q^{14};q^{14})_\infty}{(q;q)_\infty},  
\end{align}
\end{subequations}
\begin{subequations}\label{ssb}
\begin{align}
\label{s81}\sum_{k=0}^\infty \frac{q^{k(k+1)/2}}{(q;q)_k(q;q^2)_k}
&=\frac{(q,q^6,q^7;q^7)_\infty(q^5,q^9;q^{14})_\infty}{(q;q)_\infty(q;q^2)_\infty},
 \\
\label{s80}\sum_{k=0}^\infty \frac{q^{k(k+1)/2}}{(q;q)_k(q;q^2)_{k+1}}
&=\frac{(q^2,q^5,q^7;q^7)_\infty(q^3,q^{11};q^{14})_\infty}{(q;q)_\infty(q;q^2)_\infty}, \\
\label{s82}\sum_{k=0}^\infty \frac{q^{k(k+3)/2}}{(q;q)_k(q;q^2)_{k+1}}
&=\frac{(q^3,q^4,q^7;q^7)_\infty(q,q^{13};q^{14})_\infty}{(q;q)_\infty(q;q^2)_\infty}, \\
\label{s117} \sum_{k=0}^\infty\frac{q^{k^2}}{(q^4;q^4)_k(q;q^2)_k}
&=\frac{(q^3,q^{11},q^{14};q^{14})_\infty(q^8,q^{20};q^{28})_\infty}{(q;q^2)_\infty(q^4;q^4)_\infty},\\
\label{s118} \sum_{k=0}^\infty\frac{q^{k(k+2)}}{(q^4;q^4)_k(q;q^2)_k}
&=\frac{(q,q^{13},q^{14};q^{14})_\infty(q^{12},q^{16};q^{28})_\infty}{(q;q^2)_\infty(q^4;q^4)_\infty},\\
\label{s119}\sum_{k=0}^\infty\frac{q^{k(k+2)}}{(q^4;q^4)_k(q;q^2)_{k+1}}
&=\frac{(q^5,q^{9},q^{14};q^{14})_\infty(q^4,q^{24};q^{28})_\infty}{(q;q^2)_\infty(q^4;q^4)_\infty}.
\end{align}
\end{subequations}
They appear in Slater's list \cite{slater} as, respectively, number 
33, 32, 31, 61, 60, 59, 81, 80, 82, 117, 118 and 119. The identities 
\eqref{ssa} and \eqref{s81}--\eqref{s82} were all found by Rogers 
\cite{r94,r17}, whereas \eqref{s117}--\eqref{s119} are due to Slater. 
Since \eqref{s33}--\eqref{s31} were independently found by  Atle Selberg (as a high school student) \cite{sel36}, they are often called the Rogers--Selberg identities.
Although he does not state them explicitly, it appears from the final sentences of 
\cite{sel38} that Selberg also knew 
 \eqref{s61}--\eqref{s59}.
 The literature on these identities is large, and we refer the reader to  \cite{ag89,and80,and09,bailey,bress80,chu,dys43,hahn03,schultz,sel36,sel38,sills03,sills08} for a selection of relevant papers.
 
 Although one can  give shorter proofs of each individual identity, one benefit of our approach is that we obtain all six identities in each group, \eqref{ssa} and \eqref{ssb}, 
simultaneously.  Moreover, the recursions arising from the proof can be interpreted in terms of partitions and may have some independent interest.

\section{Preliminaries}\label{ps}

We will use the standard notation
$$(a;q)_k=\prod_{j=0}^{k-1}(1-aq^j), $$
where $k$ may be infinite. This is extended to negative $k$ through
$(a;q)_{-k}=1/(aq^{-k};q)_k$. We also introduce the multiplicative theta function
$$\theta(z;q)=(q,z,q/z;q)_\infty=\sum_{k=-\infty}^\infty(-1)^kq^{\frac{k(k-1)}2} z^k $$
and write
$$(a_1,\dots,a_m;q)_k=(a_1;q)_k\dotsm (a_m;q)_k, $$
$$\theta(a_1,\dots,a_m;q)=\theta(a_1;q)\dotsm\theta(a_m;q). $$

We will need the classical theta function identities
\begin{equation}\label{thd}\theta(z;q)=\theta(-qz^2;q^4)-z\theta(-q^3z^2;q^4),
\end{equation}
\begin{equation}\label{thm}
 \theta(xy,x/y;q)=\theta(-x^2,-qy^2;q^2)-\frac x y\,\theta(-qx^2,-y^2;q^2),
\end{equation}
 Euler's two $q$-exponential
summations
\begin{equation}\label{ei}\sum_{k=0}^\infty\frac{z^k}{(q;q)_\infty}=\frac 1{(z;q)_\infty},
\end{equation}
\begin{equation}\label{eib}\sum_{k=0}^\infty\frac{q^{\frac{k(k-1)}2}z^k}{(q;q)_k}=(-z;q)_\infty
\end{equation}
and the identity \cite{and86}
\begin{equation}\label{abs}\sum_{k=-\infty}^\infty\frac{(-1)^kq^{\frac{k(k-1)}2}}{(t;q)_kz^k}=\frac{\theta(1/z;q)}{(t,tz;q)_\infty},\qquad 0<|z|<\frac 1{|t|}.\end{equation}

\section{The summations \eqref{ssa}}
It will be convenient to multiply the identities  \eqref{s33}--\eqref{s31} by $(-q;q)_\infty$ and write them as
\begin{subequations}\label{abcev}
\begin{align}
\label{aev}A(q)&=(-q;q)_\infty\sum_{k=0}^\infty\frac{q^{2k^2}}{(q^4;q^4)_k(-q;q^2)_k}=\frac{1}{(q,q^2,q^5,q^6;q^7)_\infty},\\
B(q)&=(-q;q)_\infty\sum_{k=0}^\infty\frac{q^{2k(k+1)}}{(q^4;q^4)_k(-q;q^2)_k}=\frac{1}{(q,q^3,q^4,q^6;q^7)_\infty},\\
\label{cev}C(q)&=(-q;q)_\infty\sum_{k=0}^\infty\frac{q^{2k(k+1)}}{(q^4;q^4)_k(-q;q^2)_{k+1}}=\frac{1}{(q^2,q^3,q^4,q^5;q^7)_\infty}.
\end{align}
\end{subequations}
The first step in the proof is to show that the series 
 \eqref{s61}--\eqref{s59} can be expressed as an elementary factor times $A(q^2)$, $B(q^2)$ and $C(q^2)$, respectively. 
 Since this relation is  obvious for the product sides, it follows that
  \eqref{s61}--\eqref{s59} are pairwise equivalent to
  \eqref{s33}--\eqref{s31}. 
 In the following result we provide one-parameter extensions of the necessary identities.
 The equation
 \eqref{4ta} is the limit case $b=0$, $n\rightarrow\infty$ of \cite[(3.19)]{bw}
 and \eqref{4tb} arises in a similar way from an unpublished companion to  that identity \cite{w2}.
We provide a direct proof using Andrews's constant term method \cite[\S 4.2]{and86}.

 \begin{lemma}\label{qtl}
 The following quartic transformation formulas hold:
 \begin{align}\label{4ta}\sum_{k=0}^\infty\frac{q^{k(k-1)}t^k}{(q;q)_k(t;q^2)_k}&=\frac{(-qt;q^2)_\infty}{(t;q^2)_\infty}\sum_{k=0}^\infty\frac{q^{2k(2k-1)}t^{2k}}{(q^4;q^4)_k(-q t;q^2)_{2k}}, \\
 \label{4tb}\sum_{k=0}^\infty\frac{q^{k(k-1)}t^k}{(q;q)_k(qt;q^2)_k}&=\frac{(-t;q^2)_\infty}{(qt;q^2)_\infty}\sum_{k=0}^\infty\frac{q^{4k^2}t^{2k}}{(q^4;q^4)_k( -t;q^2)_{2k}}. \end{align}
 \end{lemma}

\begin{proof}
We start from the observation that
$$
\frac{\theta(1/z;q^2)}{(tz;q^2)_\infty(-tz;q)_\infty}=
\frac{\theta(1/z;q^2)}{(tz,-tz,-qtz;q^2)_\infty}
=\frac{\theta(1/z;q^2)}{(t^2z^2;q^4)_\infty(-qtz;q^2)_\infty}.$$
Applying
 \eqref{ei} and \eqref{abs}, it follows that
\begin{multline*} (t;q^2)_\infty\sum_{m=-\infty}^\infty\frac{(-1)^m q^{m(m-1)}}{(t;q^2)_mz^m}
\sum_{k=0}^\infty\frac{(-tz)^k}{(q;q)_k}\\
= (-qt;q^2)_\infty\sum_{m=-\infty}^\infty\frac{(-1)^mq^{m(m-1)}}{(-qt;q^2)_mz^m}
\sum_{k=0}^\infty\frac{(tz)^{2k}}{(q^4;q^4)_k},\qquad 0<|z|<1/|t|.
 \end{multline*}
Identifying 
the coefficient of $z^0$  gives \eqref{4ta}.
The identity \eqref{4tb} is obtained in the same way from
$$ \frac{\theta(1/z;q^2)}{(qtz;q^2)_\infty(-tz;q)_\infty}
=\frac{\theta(1/z;q^2)}{(q^2t^2z^2;q^4)_\infty(-tz;q^2)_\infty}.$$
\end{proof}

Substituting  $t=q$ in \eqref{4ta},  
$t=q^2$ in \eqref{4tb} and  $t=q^3$ in \eqref{4ta} 
gives the following identities.  As we have mentioned, they show that
\eqref{s61}--\eqref{s59} are equivalent to
  \eqref{s33}--\eqref{s31}, so all six identities will follow from \eqref{abcev}.

\begin{corollary}
We have
\begin{subequations}\label{ar}
\begin{align}
\label{aa}A(q^2)&=(q;q^2)_\infty\sum_{k=0}^\infty\frac{q^{k^2}}{(q;q)_k(q;q^2)_k},\\
B(q^2)&=(q;q^2)_\infty\sum_{k=0}^\infty\frac{q^{k(k+1)}}{(q;q)_k(q;q^2)_{k+1}}
,\\
C(q^2)&=(q;q^2)_\infty\sum_{k=0}^\infty\frac{q^{k(k+2)}}{(q;q)_k(q;q^2)_{k+1}}.
\end{align}
\end{subequations}
\end{corollary}

 To obtain functional equations for $A$, $B$, $C$ we will need the following 
 three-term quadratic transformation formula. It can be obtained from
\cite[Thm.~A\textsubscript{1}]{and66} by letting
$b,t\rightarrow 0$ and $a\rightarrow\infty$ with $at$ fixed.

\begin{lemma}
One has
\begin{multline}\label{ttqu}(-x;q)_\infty \sum_{k=0}^\infty \frac{q^{k(k-1)}y^k}{(q^2;q^2)_k(-x;q)_k}\\
=(-y;q^2)_\infty \sum_{k=0}^\infty \frac{q^{k(2k-1)}x^{2k}}{(q;q)_{2k}(-y;q^2)_k}
+(-qy;q^2)_\infty \sum_{k=0}^\infty \frac{q^{k(2k+1)}x^{2k+1}}{(q;q)_{2k+1}(-qy;q^2)_k}.\end{multline}
\end{lemma}

Although it may look complicated, it is very easy to prove \eqref{ttqu}.
 It will be useful to indicate a proof of a more general result.
 Namely, starting from the series
$$\sum_{k,m=0}^\infty\frac{q^{\frac{k(k-1)}2}p^{\frac{m(m-1)}2}r^{km}x^ky^m}{(q;q)_k(p;p)_m} $$
we can use \eqref{eib} to sum in either $k$ or $m$ and obtain
\begin{equation}\label{tbt}
\sum_{k=0}^\infty \frac{p^{\frac{k(k-1)}2}y^k}{(p;p)_k}(-r^kx;q)_\infty=
\sum_{k=0}^\infty \frac{q^{\frac{k(k-1)}2}x^k}{(q;q)_k}(-r^ky;p)_\infty. 
\end{equation}
Then, \eqref{ttqu} is simply the special case $p=q^2$, $r=q$, where we  have split the sum on the right  into terms with $k$ even and $k$ odd.

\begin{corollary}\label{fel}
The series $A$, $B$ and $C$ satisfy
\begin{subequations}\label{fe}
\begin{align}
\label{afe}A(q)&=\frac 1{(q^2;q^4)_\infty^2}\left(qC(-q^2)+A(q^8)\right)
,\\
\label{bfe} B(q)&=\frac 1{(q^2;q^4)_\infty^2}\left(A(-q^2)+qB(q^8)\right),\\
 \label{cfe} C(q)&=\frac 1{(q^2;q^4)_\infty^2}\left(B(-q^2)+q^3C(q^8)\right).
\end{align}
\end{subequations}
\end{corollary}

\begin{proof}
Substituting $(q,x,y)\mapsto (q^2,q,q^2)$ in \eqref{ttqu} gives
\begin{multline*}
(-q;q^2)_\infty\sum_{k=0}^\infty\frac{q^{2k^2}}{(q^4;q^4)_k(-q;q^2)_k}\\
={(-q^2;q^4)_\infty}
\sum_{k=0}^\infty\frac{q^{4k^2}}{(q^4;q^4)_k(q^4;q^8)_k}+
q(-q^4;q^4)_\infty\sum_{k=0}^\infty\frac{q^{4k(k+1)}}{(q^8;q^8)_k(q^2;q^4)_{k+1}}.\end{multline*}
Recognizing the three series from \eqref{aev}, \eqref{aa} and \eqref{cev}, we
can write it in the form
\eqref{afe}. The other two identities are obtained in the same way.
 \end{proof}
 
We  observe that \eqref{fe} can be viewed as a recursion for computing the power series coefficients
of $A(q)$, $B(q)$ and $C(q)$, starting from the initial values
$A(0)=B(0)=C(0)=1$. Hence, 
to complete the proof of \eqref{s33}--\eqref{s59}
we only need to verify the following fact.

\begin{lemma}
The equations \eqref{fe} hold for the right-hand sides of \eqref{abcev}.
\end{lemma}

\begin{proof}
Explicitly, the right-hand side of 
\eqref{afe} is
$$\frac 1{(q^2;q^4)_\infty^2}
 \left(\frac q{(q^4,-q^6,q^8,-q^{10};-q^{14})_\infty}
+\frac 1{(q^8,q^{16},q^{40},q^{48};q^{56})_\infty}\right). $$
By standard manipulation of infinite products, such as
$$\frac 1{(q^4;-q^{14})_\infty}=\frac{(-q^4;-q^{14})_\infty}{(q^8;q^{28})_\infty}
=\frac{(-q^4,q^{18};q^{28})_\infty}{(q^8;q^{28})_\infty}, $$
$$\frac 1{(q^8;q^{56})_\infty}=\frac{(q^{36};q^{56})_\infty}{(q^8;q^{28})_\infty}
= \frac{(q^{18},-q^{18};q^{28})_\infty}{(q^8;q^{28})_\infty},$$
this can be written
\begin{multline*}\frac{(q^6,q^{10},q^{18},q^{22};q^{28})}{(q^2;q^4)_\infty^2(q^8,q^{12},q^{16},q^{20};q^{28})_\infty}\\
\times\left(q(-q^4,-q^{8},-q^{20},-q^{24};q^{28})_\infty+(-q^6,-q^{10},-q^{18},-q^{22};q^{28})_\infty\right).
\end{multline*}
By \eqref{thm} with $(q,x,y)\mapsto (q^{14},q^3, -q^2)$, the final factor can be expressed as
$$\frac {(q^{14};q^{14})_\infty^2}{(q^{28};q^{28})_\infty^2}
(-q,-q^5,-q^9,-q^{13};q^{14})_\infty
$$
After some further  manipulation, this leads to the product formula for $A(q)$.
The proofs of \eqref{bfe} and \eqref{cfe} are similar, using the cases
$(q,x,y)\mapsto(q^{14},q^2,-q)$ and $(q^{14},q^4,-q)$ of \eqref{thm}, respectively.
\end{proof}

We conclude this section by noting that the equations \eqref{fe} have a combinatorial 
interpretation. If we write
$$A(q)=\sum_{k=0}^\infty a_kq^k, $$
then it is clear from \eqref{aev} that $a_k$ is the number of partitions of $k$ into
parts congruent to $\pm 1$ or $\pm 2$ mod $7$. Analogous statements hold for $B$ and $C$. Moreover, the coefficients in
$$\frac 1{(q;q^2)_\infty^2}=(-q;q)_\infty^2=\sum_{k=0}^\infty \phi_k q^k $$
enumerate  $2$-colored partitions into odd parts or, alternatively, into distinct parts. This leads to the following reformulation of \eqref{fe}.

\begin{corollary}\label{cic}
Let $a_k$ be the number of partitions of $k$ into parts congruent to
 $\pm 1$ or $\pm 2$ mod $7$, $b_k$ the number of partitions of $k$ into parts congruent to
 $\pm 1$ or $\pm 3$ mod $7$, $c_k$ the number of partitions of $k$ into parts congruent to
 $\pm 2$ or $\pm 3$ mod $7$ and $\phi_k$  the 
 number of $2$-colored partitions of $k$ into odd parts.
 Then,
 \begin{align*}
 a_{2k}&=\sum_{m=0}^{\lfloor k/4\rfloor}\phi_{k-4m} a_m, & 
 a_{2k+1}&=\sum_{m=0}^k(-1)^m\phi_{k-m} c_m,\\
 b_{2k}&=\sum_{m=0}^k(-1)^m\phi_{k-m} a_m, & b_{2k+1}&=\sum_{m=0}^{\lfloor k/4\rfloor}\phi_{k-4m} b_m,\\
 c_{2k}&=\sum_{m=0}^k(-1)^m\phi_{k-m} b_m, & c_{2k+1}&=\sum_{m=0}^{\lfloor (k-1)/4\rfloor}\phi_{k-4m-1} c_m.
 \end{align*}
\end{corollary}


\section{The summations \eqref{ssb}}

The proof of \eqref{ssb} is very similar to that of \eqref{ssa}.
We first write \eqref{s81}--\eqref{s82} as
\begin{subequations}\label{def}
\begin{align}\label{d}D(q)&=(q;q^2)_\infty
\sum_{k=0}^\infty \frac{q^{k(k+1)/2}}{(q;q)_k(q;q^2)_k}
=\frac{1}{(q^2,q^3,q^4,q^{10},q^{11},q^{12};q^{14})_\infty},
 \\
E(q)&=(q;q^2)_\infty\sum_{k=0}^\infty \frac{q^{k(k+1)/2}}{(q;q)_k(q;q^2)_{k+1}}
=\frac{1}{(q,q^4,q^6,q^{8},q^{10},q^{13};q^{14})_\infty}, \\
F(q)&=(q;q^2)_\infty\sum_{k=0}^\infty \frac{q^{k(k+3)/2}}{(q;q)_k(q;q^2)_{k+1}}
=\frac{1}{(q^2,q^5,q^6,q^{8},q^{9},q^{12};q^{14})_\infty}.
\end{align}
\end{subequations}

To relate these series  to \eqref{s117}--\eqref{s119}
we need the  case  $p=r=q^2$ of \eqref{tbt}, that is,
\begin{equation}\label{tpt}
(-x;q)_\infty\sum_{k=0}^\infty \frac{q^{k(k-1)}y^k}{(q^2;q^2)_k(-x;q)_{2k}}=
(-y;q^2)_\infty \sum_{k=0}^\infty \frac{q^{\frac{k(k-1)}2}x^k}{(q;q)_k(-y;q^2)_k}.
\end{equation}
This can be obtained from
\cite[Thm.~A\textsubscript{3}]{and66} by letting
$b,t\rightarrow 0$ and $a\rightarrow\infty$ with $at$ fixed.
Specializing the parameters in \eqref{tpt}, we find that 
\begin{subequations}\label{defa}
\begin{align}
\label{da}D(q)&=\frac 1{(q;q^2)_\infty}\sum_{k=0}^\infty\frac{(-1)^kq^{k^2}}{(q^4;q^4)_k(-q;q^2)_k},\\
E(q)&=\frac 1{(q;q^2)_\infty}\sum_{k=0}^\infty\frac{(-1)^kq^{k(k+2)}}{(q^4;q^4)_k(-q;q^2)_k},\\
\label{fa}F(q)&=\frac 1{(q;q^2)_\infty}\sum_{k=0}^\infty\frac{(-1)^kq^{k(k+2)}}{(q^4;q^4)_k(-q;q^2)_{k+1}},
\end{align}
\end{subequations}
where we recognize the series in 
\eqref{s117}--\eqref{s119} with $q\mapsto -q$. It is easy to check that the same relations hold between the  right-hand sides. 
Thus, the identities  \eqref{s117}--\eqref{s119} are
equivalent to  \eqref{s81}--\eqref{s82}.

To obtain
functional equations for \eqref{def} we use the following 
three-term quartic transformation formulas, which we have not found in the literature.

\begin{lemma}\label{ttql}
The following identities hold:
\begin{subequations}\label{ttq}
\begin{align}
\nonumber(t;q)_\infty\sum_{k=0}^\infty\frac{q^{\frac{k(k-1)}2}t^k}{(q^2;q^2)_k(t;q)_k}
\label{ttqa}&=(-qt^2;q^4)_\infty\sum_{k=0}^\infty\frac{q^{k(2k+1)}t^{2k}}{(q^2;q^2)_{2k}(-qt^2;q^4)_k}\\
&\quad-(-q^3t^2;q^4)_\infty\sum_{k=0}^\infty\frac{q^{(k+1)(2k+1)}t^{2k+1}}{(q^2;q^2)_{2k+1}(-q^3t^2;q^4)_k},\\
\nonumber(t;q)_\infty\sum_{k=0}^\infty\frac{q^{\frac{k(k+1)}2}t^k}{(q^2;q^2)_k(t;q)_k}
&=(-q^3t^2;q^4)_\infty\sum_{k=0}^\infty\frac{q^{k(2k-1)}t^{2k}}{(q^2;q^2)_{2k}(-q^3t^2;q^4)_k}\\
\label{ttqb}&\quad-(-q^5t^2;q^4)_\infty\sum_{k=0}^\infty\frac{q^{k(2k+1)}t^{2k+1}}{(q^2;q^2)_{2k+1}(-q^5t^2;q^4)_k}.
\end{align}
\end{subequations}
\end{lemma} 

\begin{proof}
Using \eqref{ei} and \eqref{abs},
we can write the left-hand side of \eqref{ttqa} as
$$(t;q)_\infty\ct\sum_{k=-\infty}^\infty\frac{(-1)^kq^{\frac{k(k-1)}{2}}}{(t;q)_kz^k}
\sum_{m=0}^\infty \frac{(-tz)^m}{(q^2;q^2)_m}
=\ct\frac{\theta(1/z;q)}{(tz;q)_\infty(-tz;q^2)_\infty}, $$
where $\ct$ denote the constant term in the Laurent expansion for $z$ near $0$.
By
\eqref{thd}, the function on the right is
\begin{multline*}\frac{\theta(-q/z^2;q^4)-z^{-1}\theta(-q^3/z^2;q^4)}{(t^2z^2;q^4)(qtz;q^2)_\infty}
= (-qt^2;q^4)_\infty\sum_{k=-\infty}^{\infty}\frac{q^{k(2k-1)}}{(-qt^2;q^4)_k z^{2k}}\sum_{m=0}^\infty \frac{(qtz)^m}{(q^2;q^2)_m}\\
- (-q^3t^2;q^4)_\infty\sum_{k=-\infty}^{\infty}\frac{q^{k(2k+1)}}{(-q^3t^2;q^4)_k z^{2k+1}}\sum_{m=0}^\infty \frac{(qtz)^m}{(q^2;q^2)_m}. \end{multline*}
Picking out the constant term we arrive at the right-hand side of \eqref{ttqa}.
The identity \eqref{ttqb} is obtained in the same way, starting from
$$\ct\frac{\theta(1/z;q)}{(tz;q)_\infty(-qtz;q^2)_\infty}.$$
\end{proof}

Replacing $(q,t)\mapsto(q^2,-q)$ in  \eqref{ttqa} 
gives
\begin{multline*}
(-q;q^2)_\infty\sum_{k=0}^\infty\frac{(-1)^kq^{k^2}}{(q^2,-q,-q^2;q^2)_k}
=(-q^4;q^8)_\infty\sum_{k=0}^\infty\frac{q^{4k(k+1)}}{(q^4,-q^4,q^8;q^8)_k}\\
+\frac{q^3}{1-q^4}(-q^8;q^8)_\infty\sum_{k=0}^\infty\frac{q^{4k(k+2)}}{(q^8,-q^8,q^{12};q^8)_k},
\end{multline*}
where we  recognize the three series from \eqref{da}, \eqref{d} and \eqref{fa}, respectively. Considering also the cases 
 $(q,t)\mapsto(q^2,-q)$ of \eqref{ttqb} and $(q,t)\mapsto(q^2,-q^3)$ of \eqref{ttqa} leads to the following result.

  \begin{lemma}
  The series \eqref{def} satisfy the system of functional equations
  \begin{subequations}\label{deff}
  \begin{align}
\label{df}  D(q)&=\frac{(q^8;q^8)_\infty}{(q^2;q^2)_\infty} \left(D(q^8)+q^3 F(-q^4)\right),\\
 \label{ef} E(q)&= \frac{(q^8;q^8)_\infty}{(q^2;q^2)_\infty} \left( q E(q^8)+D(-q^4)\right),\\
 \label{ff}  F(q)&=\frac{(q^8;q^8)_\infty}{(q^2;q^2)_\infty} \left( q^5 F(q^8)+ E(-q^4)\right).
  \end{align}
  \end{subequations}
  \end{lemma}

As before, the equations \eqref{deff} determine $D$, $E$ and $F$ as power series, up to the normalizing condition $D(0)=E(0)=F(0)=1$. Thus, all that remains to prove \eqref{ssb} is the following fact.

\begin{lemma} \label{dpfl}
The identities \eqref{deff} hold for the right-hand sides of \eqref{def}.
\end{lemma}

\begin{proof}
Using now  $D$, $E$ and $F$ to denote the product expressions, we 
observe that
\begin{align*}
D(q)&=\frac{(q,q^6,q^7;q^7)_\infty(q^5,q^9;q^{14})_{\infty}}{(q;q)_\infty}=\frac{(q,q^6,q^7,q^{5/2},-q^{5/2},q^{9/2},-q^{9/2};q^7)_\infty}{(q;q)_\infty}\\
&=\frac {(q,q^{5/2};q^{7/2})_\infty(q^7,-q^{5/2},-q^{9/2};q^7)_\infty}{(q;q)_\infty},\\
F(-q)&=\frac 1{(q^2,-q^5,q^6,q^8,-q^9,q^{12};q^{14})_\infty}
=\frac {(-q^2,-q^{12};q^{14})_\infty}{(q^2,q^6,q^8,q^{12};q^{14})_\infty(-q^2,-q^5;q^7)_\infty}\\
&=\frac {(q^2,q^5;q^7)_\infty(-q^2,-q^{12};q^{14})_\infty}{(q^2,q^4,q^6,q^8,q^{10},q^{12};q^{14})_\infty}=\frac {(q^2,q^5;q^7)_\infty(-q^2,-q^{12},q^{14};q^{14})_\infty}{(q^2;q^{2})_\infty}.
\end{align*}
It follows that the right-hand side of \eqref{df} is
$$\frac{(q^8,q^{20};q^{28})_\infty(q^{56};q^{56})_\infty}{(q^2;q^2)_\infty}\left(
(-q^{20},-q^{36};q^{56})_\infty+q^3(-q^8,-q^{48};q^{56})_\infty
\right). $$
By \eqref{thd} with $(q,z)\mapsto(q^{14},-q^3)$, this equals
\begin{align*}\frac{(q^8,q^{20};q^{28})_\infty(-q^3,-q^{11},q^{14};q^{14})_\infty}{(q^2;q^2)_\infty} 
&=\frac{(q^6,q^8,q^{14};q^{14})_\infty}{(q^2;q^2)_\infty(q^3,q^{11};q^{14})_\infty}\\
&=\frac 1{(q^2,q^3,q^4,q^{10},q^{11},q^{12};q^{14})_\infty}.
\end{align*}
This gives \eqref{df}. In the same way, \eqref{ef} and \eqref{ff}
follow from  \eqref{thd} with $(q,z)\mapsto(q^{14},-q)$ and $(q^{14},-q^5)$, respectively.
\end{proof}

Finally, we give a combinatorial reformulation of Lemma \ref{dpfl}.

\begin{corollary}
Let $d_k$ denote the number of partitions of $k$ into parts congruent to $\pm 3$ mod $7$ or
 $\pm 2$ mod $14$, $e_k$ the number of partitions of $k$ into parts congruent to $\pm 1$ mod $7$ or $\pm 4$ mod $14$, $f_k$ the number of partitions of $k$ into parts congruent to $\pm 2$ mod $7$ or $\pm 6$ mod $14$ and 
 $\psi_k$ the number of partitions of $k$ into parts not divisible by $4$. Then,
 \begin{align*}
 d_{2k}&=\sum_{m=0}^{\lfloor k/4\rfloor}\psi_{k-4m} d_m, & 
 d_{2k+1}&=\sum_{m=0}^{\lfloor(k-1)/2\rfloor}(-1)^m\psi_{k-2m-1} f_m,\\
 e_{2k}&=\sum_{m=0}^{\lfloor k/2\rfloor}(-1)^m\psi_{k-2m} d_m, & e_{2k+1}&=\sum_{m=0}^{\lfloor k/4\rfloor}\psi_{k-4m} e_m,\\
 f_{2k}&=\sum_{m=0}^{\lfloor k/2\rfloor}(-1)^m\psi_{k-2m} e_m, & f_{2k+1}&=\sum_{m=0}^{\lfloor (k-2)/4\rfloor}\psi_{k-4m-2} f_m.
 \end{align*}
\end{corollary}

\end{document}